\documentclass[12pt]{article}
\usepackage{amssymb}
\usepackage{latexsym}
\usepackage{pstcol, pst-plot}

\newcommand{\ben}{\begin{enumerate}}
\newcommand{\een}{\end{enumerate}}
\newcommand{\bde}{\begin{defn}}
\newcommand{\ede}{\end{defn}}
\newcommand{\bex}{\begin{exa}}
\newcommand{\eex}{\end{exa}}
\newcommand{\barr}{\begin{array}}
\newcommand{\earr}{\end{array}}
\newcommand{\btab}{\begin{tabular}}
\newcommand{\etab}{\end{tabular}}
\newcommand{\beq}{\begin{equation}}
\newcommand{\eeq}{\end{equation}}
\newcommand{\bea}{\begin{eqnarray*}}
\newcommand{\eea}{\end{eqnarray*}}
\newcommand{\bce}{\begin{center}}
\newcommand{\ece}{\end{center}}

\newcommand{\bib}{thebibliography}

\newcommand{\qed}{\mbox{$\Diamond$}\vspace{\baselineskip}}
\newcommand{\qedthm}{\quad\mbox{$\Diamond$}}
\newtheorem{theorem}{Theorem}[section]
\newtheorem{proposition}[theorem]{Proposition}
\newtheorem{lemma}[theorem]{Lemma}
\newtheorem{definition}[theorem]{Definition}
\newtheorem{corollary}[theorem]{Corollary}

\newtheorem{conjecture}[theorem]{Conjecture}
\newenvironment{proof}{\noindent {\bf Proof:}}{{\qed}}
\newcommand{\bfi}{\begin{figure} \begin{center}}
\newcommand{\efi}{\end{center} \end{figure}}

\newcommand{\ble}{\begin{lemma}}
\newcommand{\ele}{\end{lemma}}
\newcommand{\bth}{\begin{theorem}}
\renewcommand{\eth}{\end{theorem}}
\newcommand{\bpr}{\begin{proposition}}
\newcommand{\epr}{\end{proposition}}
\newcommand{\bco}{\begin{corollary}}
\newcommand{\eco}{\end{corollary}}
\newcommand{\bcon}{\begin{conjecture}}
\newcommand{\econ}{\end{conjecture}}

\newcommand{\al}{\alpha}
\newcommand{\be}{\beta}
\newcommand{\ga}{\gamma}
\newcommand{\de}{\delta}
\newcommand{\ep}{\epsilon}
\newcommand{\ree}[1]{(\ref{#1})}

\makeatletter
\newfont{\footsc}{cmcsc10 at 8truept}
\newfont{\footbf}{cmbx10 at 8truept}
\newfont{\footrm}{cmr10 at 10truept}
\newcommand*{\ppat}{\nolinebreak{1(\ell+1)\ell\dots 2}}
\newcommand*{\permpat}{q_\ell}

\begin{document}
\pagestyle{plain}
\title{Bounding quantities related to the packing density of $\ppat$}
\author{Martin Hildebrand\\
\small Department of Mathematics and Statistics\\[-0.8ex]
\small State University of New York - University at Albany, Albany, NY 
12222\\[-0.8ex]
\small \texttt{martinhi@math.albany.edu}\\[1.6ex]
Bruce E. Sagan\\
\small Department of Mathematics\\[-0.8ex]
\small Michigan State University, East Lansing, MI 48824\\[-0.8ex]
\small \texttt{sagan@math.msu.edu}\\[1.6ex]
Vincent R. Vatter\thanks{Partially supported by an NSF VIGRE grant to 
the Rutgers University Department of Mathematics.}\\
\small Department of Mathematics\\[-0.8ex]
\small Rutgers University, Piscataway, NJ 08854\\[-0.8ex]
\small \texttt{vatter@math.rutgers.edu}}

\date{\today}
\maketitle

\begin{abstract}
We bound several quantities related to the packing density of the patterns $\ppat$.  These bounds sharpen results of B\'ona, Sagan, and Vatter and give a new proof of the packing density of these patterns, originally computed by Stromquist in the case $\ell=2$ and by Price for larger $\ell$.  We end with comments and conjectures.
\end{abstract}

\pagestyle{plain}

\section{Introduction}

We say two sequences $p, q$ of length $n$ are of the same {\it type\/} 
if $p(i)<p(j)$ if and only if $q(i)<q(j)$ for all $i,j\in[n]$, that is, 
if $p$ and $q$ have the same pairwise comparisons.  For an 
$n$-permutation $p$ and an $\ell$-permutation $q$ we let $c_q(p)$ denote 
the number of $l$-subsequences of type $q$ in $p$, and we say that $p$ 
contains $c_q(p)$ copies of the {\it pattern\/} $q$.  For example, 41523
contains exactly  
two 132-patterns, namely 152 and 153, so $c_{132}(41523)=2$.

We say that an $n$-permutation $p$ is {\it $q$-optimal\/} if there is no 
$n$-permutation with more  copies of $q$ than $p$, and let
$$
\mbox{$M_{n,q}=c_q(p)$ for a $q$-optimal $p$}.
$$
Since there are a total of ${n\choose \ell}$ $\ell$-subsequences in any 
$n$-permutation, we always have $0\le M_{n,q}\le {n\choose |q|}$.  The {\it 
packing density\/} of a permutation $q$ is defined as
$$
\delta(q)=\lim_{n\rightarrow\infty} \frac{M_{n,q}}{{n\choose |q|}}.
$$
This limit exists because of the following theorem.  An unpublished
proof was given by Galvin and reproduced in Price's
thesis.  One can also find the demonstration in a paper
of Albert, Atkinson, Handley, Holton, and Stromquist.
\bth\cite{aahhs:opdop, price} 
\label{albert}
The ratio  $M_{n,q}/{n\choose |q|}$ is weakly decreasing.\qedthm
\eth

Stromquist~\cite{str:plp} computed the packing density of $132$.
Using similar techniques, Price computed the packing density of 
the patterns $q_\ell=\ppat$ for all
$\ell\ge2$. 

\begin{theorem}\label{price}\cite{price}
The packing density of $q_\ell$ is 
\beq
\label{beta}
\beta=\ell\alpha(1-\alpha)^{\ell-1}
\eeq
where $\alpha$ is the unique root of 
\beq
\label{f_ell}
f_\ell(x)=\ell x^{\ell+1}-(\ell+1)x+1
\eeq
in the interval $(0,1)$.\qedthm
\end{theorem}

Since $f_\ell(1/(\ell+1))>0$ and $f_\ell(1/\ell)<0$, we have
\beq\label{alpha-bounds}
\frac{1}{\ell+1}<\alpha<\frac{1}{\ell}.
\eeq
The chart below shows approximate 
values of $\alpha$ and $\beta$ for small $\ell$.
$$
\barr{c|c|c}
\ell	&\alpha		&\beta\\
\hline \hline
2	&0.366	&0.464\rule{0pt}{20pt}\\[5pt]
\hline
3	&0.253	&0.424\rule{0pt}{20pt}\\[5pt]
\hline
4	&0.200	&0.410\rule{0pt}{20pt}\\[5pt]
\hline
5	&0.167	&0.402\rule{0pt}{20pt}\\[5pt]
\earr
$$

For the rest of the paper, we abbreviate $M_{n,q_\ell}$ to $M_n$.  Price 
proved Theorem~\ref{price} by showing that
$$
\frac{M_n}{{n\choose \ell+1}}=\beta + O\left(\frac{\log n}{n}\right).
$$
We will reprove Theorem~\ref{price} by giving precise bounds on 
$M_n$.
\begin{theorem}\label{MAIN}
For all $n\ge\ell\ge 2$,
$$
\beta\frac{(n-\ell)^{\ell+1}}{(\ell+1)!}\le M_n\le 
\beta\frac{(n+\de_{2,\ell})^{\ell+1}}{(\ell+1)!},
$$
where $\de_{2,\ell}$ is the Kronecker delta (and not to be confused
with the packing density $\beta=\de(q_\ell)$).
\end{theorem}

Note that Theorem~\ref{price} follows immediately from the theorem
just stated by merely dividing all sides by ${n\choose \ell+1}$ and
taking $n\rightarrow\infty$.
Also note that the lower bound follows from Price's calculation of the packing density of $q_\ell$ and the fact that
$M_n/{n\choose \ell+1}$ is  
decreasing, but we will provide another demonstration in order to give a new proof of Theorem~\ref{price}.  We will also have other uses for the intermediate 
results needed to prove both bounds.  

The rest of this paper is structured as follows.  In the next section
we give some preliminary definitions and previous results which
will be needed for our bounds.  Section~\ref{M_n} is devoted to
proofs of bounds involving $M_n$.   In the section following that,
we provide bounds for a related quantity.
Often our upper bound proofs from these sections will not work when
$\ell=2$, so Section~\ref{ell=2} is devoted to a discussion
of that case.  Finally, we end with a section of comments and conjectures.

\section{Definitions and previous results}

We say that a permutation is {\it layered\/} if it is the concatenation 
of subwords (the {\em layers}) where the entries decrease within each 
layer, and increase between the layers.  For example, $321548769$ is a 
layered permutation with layers $321,54,876$, and 9.  The only 
permutations for which the packing density has been computed are layered 
or equivalent to layered permutations under one of the routine 
symmetries.  The following theorem of Stromquist is crucial for 
computing these densities.  Its proof may also be found in Price's 
thesis~\cite{price}, and a generalization is proved in 
\cite{aahhs:opdop}.  B\'ona, Sagan, and Vatter~\cite{bon:iz} proved a
similar result for $n$-permutations with $M_n-1$ copies of 
$q_\ell$, for any $l\ge 2$.

\begin{theorem}\cite{str:plp}\label{strom-layered}
For all layered permutations $q$ and positive integers $n$, there is a 
layered $q$-optimal $n$-permutation.\qedthm
\end{theorem}

Layered $\permpat$-optimal permutations have the following easily 
established recursive structure.

\begin{proposition}\cite{bon:iz}\label{fund}
Let $p$ be a layered $\permpat$-optimal $n$-permutation whose last layer 
is of length $m$. Then the leftmost $k=n-m$ elements of $p$ form a 
$\permpat$-optimal $k$-permutation.\qedthm
\end{proposition}

The previous proposition implies that
\beq
\label{Mn}
M_n=\max_{1\le k<n}\left(M_k+k{m \choose \ell}\right).
\eeq

The value of $k$ that maximizes the right-hand side of (\ref{Mn}) will 
be very important throughout this paper, so we give it a notation as follows.

\begin{definition} For any positive integer $n>\ell$, let $k_n$ denote 
the positive integer for which $M_k+ k {m \choose \ell}$ is maximal. If 
there are several integers with this property, let $k_n$ be the largest 
among them.
\end{definition}

Once we have found the packing density of $q_\ell$ (Theorem~\ref{price}), it is not hard to find the asymptotic behavior of $k_n$.

\begin{corollary}\cite{price}
\label{price-kn}
The limit of $k_n/n$ is $\alpha$.\qedthm
\end{corollary}

We will sharpen this result considerably in Section~\ref{k_n}.  We
will also need some information about $\be$.
First are a couple of extremal expressions for $\be$.

\begin{lemma}
\label{whatisbeta}
The quantity $\beta$ satisfies
\bea
\beta	&=&\max_{0\le\gamma\le 1}
	\frac{(\ell+1)\ga(1-\ga)^{\ell-1}}{1+\ga+\ga^2+\dots+\ga^\ell}\\
	&=&\min_{0\le\gamma\le \al}
	\be\ga^\ell+(1-\ga)^\ell.
\eea
In fact $\be=\be\al^\ell+(1-\al)^\ell$.
\end{lemma}
\begin{proof}
The maximum expression for $\be$ was
given by Price~\cite{price} in his proof of Theorem~\ref{price}.  

After rearranging
terms and plugging in the definition of $\be$, proving the last
equation is equivalent to showing that 
$$
\ell\al(1-\al)^{\ell-1}(\al^\ell-1)+(1-\al)^\ell=0.
$$
Cancelling out $(1-\al)^{\ell-1}$ leaves the defining equation for
$\al$ and thus proves the result.

Now to obtain the minimum expression, it suffices to show that 
$\be\ga^\ell+(1-\ga)^\ell$ is an decreasing function of $\ga$ on the
interval $[0,\al]$.  It is an easy exercise in calculus to show that,
in fact, it is decreasing on $[0,1/\ell]$.  So by~\ree{alpha-bounds}
we are done. 
\end{proof} 

In addition, we will need some upper bounds for $\be$.
\ble
\label{whatisbeta2}
For all $\ell\ge2$ we have
$$
\be\le\left(1-\frac{1}{\ell}\right)^{\ell-1}\le\frac{1}{2}. 
$$
\ele
\begin{proof}
For the first inequality, consider the function 
\beq
\label{f(x)}
f(x)=\ell x(1-x)^{\ell-1}.
\eeq
Clearly  $f(\al)=\be$.  Furthermore,
elementary calculus shows that $f(x)$ is an increasing function on
the interval $[0,1/\ell]$, which contains $\al$ by~\ree{alpha-bounds}.
So $f(\al)\le f(1/\ell)$ and we are done with the first bound.
For the second inequality we use the usual bounds for alternating
series to give
\bea
\left(1-\frac{1}{\ell}\right)^{\ell-1}
	&\le& 1-\frac{\ell-1}{\ell}+\frac{(\ell-1)(\ell-2)}{2\ell^2}\\
	&=&  \frac{1}{2}-\frac{1}{2\ell}+\frac{1}{\ell^2}\\[5pt]
	&\le&\frac{1}{2}
\eea
when $\ell\ge2$.
\end{proof}

B\'ona et al.\ gave crude bounds on $k_n$.
\bpr\cite{bon:iz}
\label{kn-crude}
For $n>\ell$ we have
$$
\frac{n-\ell}{\ell+1}\le k_n<\frac{n}{l}.\qedthm
$$
\epr

They also found that the sequence $k_n$ is ``continuous'' in the following
sense. 
\begin{theorem}[Continuity Theorem]\cite{bon:iz}\label{cont}
The sequence $(k_n)_{n>\ell}$ diverges to infinity and satisfies
$$
k_{n-1}\le k_n\le k_{n-1}+1
$$
for all $n> l+1$.\qedthm
\end{theorem}

The Continuity Theorem will be very useful for us because it shows that 
there are only two possibilities for $k_{n-1}$: either $k_n$ or $k_n-1$.

Let $c_{n,i}$ denote the number of copies of $q_\ell$ in an 
$n$-permutation whose last layer is of length $n-i$ and whose leftmost 
$i$ elements form a $\permpat$-optimal $i$-permutation.  So for $1\le i<n$,
\beq
\label{c_ni}
c_{n,i}=M_i+ i{n-i \choose \ell}.
\eeq
As in~\cite{bon:iz}, the sequences $(M_n)_{n\ge 1}$ and 
$(c_{n,i})_{i=1}^{n-1}$ will arise repeatedly, so we need to recall some 
results about them.
We will frequently consider the difference $c_{n,i}-c_{n,i-1}$, so let 
us simplify it now
\beq\label{cnidiff}
c_{n,i}-c_{n,i-1}=M_i-M_{i-1}+\frac{n-(\ell+1)i+1}{\ell}{n-
i\choose\ell-1}.
\eeq
We will also need the following result about differences of the $M_n$.
\begin{lemma}\cite{bon:iz}\label{DD}
For all $n\ge 0$ we have 
$$
0\le (M_{n+2}-M_{n+1})-(M_{n+1}-M_n)\le {n\choose \ell-1}.\qedthm 
$$
\end{lemma}

\begin{figure}[t]
\begin{center}
\psset{xunit=0.19in, yunit=0.0015in}
\begin{pspicture}(-2,0)(30,2500)
\psaxes[axesstyle=frame, dy=500\psyunit, Dy=500, dx=5\psxunit, Dx=5, tickstyle=bottom](0,0)(30,2500)
\psset{linewidth=0.02in}
\pscircle*(1, 406){0.05in}
\psline(1, 406)(2, 756)
\pscircle*(2, 756){0.05in}
\psline(2, 756)(3, 1053)
\pscircle*(3, 1053){0.05in}
\psline(3, 1053)(4, 1303)
\pscircle*(4, 1303){0.05in}
\psline(4, 1303)(5, 1506)
\pscircle*(5, 1506){0.05in}
\psline(5, 1506)(6, 1668)
\pscircle*(6, 1668){0.05in}
\psline(6, 1668)(7, 1791)
\pscircle*(7, 1791){0.05in}
\psline(7, 1791)(8, 1878)
\pscircle*(8, 1878){0.05in}
\psline(8, 1878)(9, 1935)
\pscircle*(9, 1935){0.05in}
\psline(9, 1935)(10, 1963)
\pscircle*(10, 1963){0.05in}
\psline(10, 1963)(11, 1968)
\pscircle*(11, 1968){0.05in}
\psline(11, 1968)(12, 1951)
\pscircle*(12, 1951){0.05in}
\psline(12, 1951)(13, 1915)
\pscircle*(13, 1915){0.05in}
\psline(13, 1915)(14, 1866)
\pscircle*(14, 1866){0.05in}
\psline(14, 1866)(15, 1806)
\pscircle*(15, 1806){0.05in}
\psline(15, 1806)(16, 1738)
\pscircle*(16, 1738){0.05in}
\psline(16, 1738)(17, 1668)
\pscircle*(17, 1668){0.05in}
\psline(17, 1668)(18, 1596)
\pscircle*(18, 1596){0.05in}
\psline(18, 1596)(19, 1527)
\pscircle*(19, 1527){0.05in}
\psline(19, 1527)(20, 1466)
\pscircle*(20, 1466){0.05in}
\psline(20, 1466)(21, 1413)
\pscircle*(21, 1413){0.05in}
\psline(21, 1413)(22, 1374)
\pscircle*(22, 1374){0.05in}
\psline(22, 1374)(23, 1353)
\pscircle*(23, 1353){0.05in}
\psline(23, 1353)(24, 1350)
\pscircle*(24, 1350){0.05in}
\psline(24, 1350)(25, 1375)
\pscircle*(25, 1375){0.05in}
\psline(25, 1375)(26, 1425)
\pscircle*(26, 1425){0.05in}
\psline(26, 1425)(27, 1504)
\pscircle*(27, 1504){0.05in}
\psline(27, 1504)(28, 1621)
\pscircle*(28, 1621){0.05in}
\psline(28, 1621)(29, 1773)
\pscircle*(29, 1773){0.05in}
\end{pspicture}
\end{center}
\caption{A plot of the sequence $(c_{30,i})_{i=1}^{29}$ when $\ell=2$.  The sequence begins at $c_{30,1}=406$ and rises to $c_{30,11}=M_{30}=1968$, so $k_{30}=11$.  It then falls until $c_{30,24}=1350$ before changing direction one last time.}\label{bimodalfig}
\end{figure}
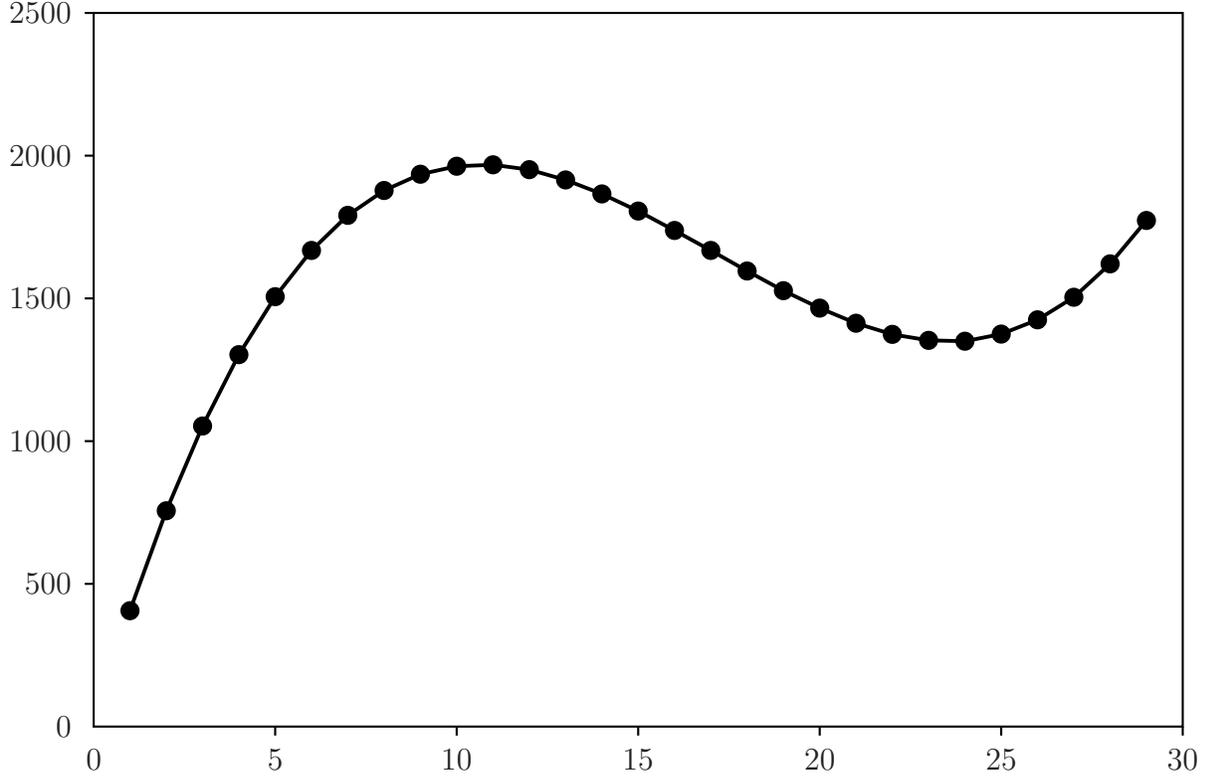

To conclude our recap of results from \cite{bon:iz}, we state the 
Bimodal Theorem.  It plays a crucial role in the arguments both in
that paper and in this one.
\begin{theorem}[Bimodal Theorem]\cite{bon:iz}\label{bimodal}
For each positive integer $n>\ell$ there is some integer $j>n/\ell$ 
(depending, of course, on $n$) so that:
\ben
\item[(i)] $c_{n,i-1} \le c_{n,i}$ if $i \le k_n$,
\item[(ii)] $c_{n,i-1} > c_{n,i}$ if $k_n < i \le j$,
\item[(iii)] $c_{n,i-1} \le c_{n,i}$ if $j<i<n$.\qedthm
\een
\end{theorem}

Figure~\ref{bimodalfig} illustrates the phenomenon described by the Bimodal Theorem.

\section{Bounds on $M_n$}
\label{M_n}

For all $k\ge 1$, let $n_k$ denote the least integer $n\ge \ell+1$ such 
that $k_n=k$.  As a trivial example, $n_1=\ell+1$.  In general, we 
always have the following upper bound on $n_k$.

\begin{proposition}\label{nk-ub}
For all $k\ge 2$, we have
$$
n_k\le(\ell+1)k-1.
$$
\end{proposition}
\begin{proof}
Substituting $n=(\ell+1)k-1$ and $i=k$ reduces (\ref{cnidiff}) to
$c_{n,k}-c_{n,k-1}=M_k-M_{k-1}\ge0$.  By the Continuity Theorem, it
suffices to show that $k_n\ge k$.  Let $j$ be as in the Bimodal 
Theorem.  By that theorem we know that $\{c_{n,i}\}_{i=1}^{n-1}$ is bimodal 
with three sections $\{c_{n,i}\}_{i=1}^{k_n}$, 
$\{c_{n,i}\}_{i=k_n}^{j}$, and $\{c_{n,i}\}_{i=j}^{n-1}$, where the first 
and last sections are weakly increasing, while the second section is 
strictly decreasing.  Therefore we must have either $k\le k_n$, as 
desired, or $k>j$.  However, $j>n/\ell=((\ell+1)k-1)/\ell>k$ for $k\ge2$, so the 
latter possibility cannot occur, finishing the proof.
\end{proof}

In the next lemma we compute $n_k$ for all sufficiently small $k$.

\begin{lemma}\label{basecases}
For all $2\le k\le \ell+1$ we have
$$
n_k=(\ell+1)k-1.
$$
\end{lemma}
\begin{proof}
Fix $k$ between $2$ and $\ell+1$.  
Consider first the case when $2\le k\le \ell$.  
Then $M_k=M_{k-1}=0$ and, by the
Continuity Theorem, we have $n_k-k>\ell-1$ so 
${n_k-k\choose \ell-1}>0$.  We use~(\ref{cnidiff}) to get 
$$
0\le c_{n_k,k}-c_{n_k,k-1}=\frac{n_k-(\ell+1)k+1}{\ell}{n_k-k\choose 
\ell-1},
$$
which yields $n_k\ge (\ell+1)k-1$.  The inequality in the other
direction is given to us by  Proposition~\ref{nk-ub}, finishing this case.

Now consider $k=\ell+1$.  By  Proposition~\ref{nk-ub} again, it
suffices to show that $n_{\ell+1}\ge (\ell+1)^2-1$.  Using (\ref{cnidiff}) again we get  
$$
c_{(\ell+1)^2-2,\ell+1}-c_{(\ell+1)^2-2,\ell}
=1-\frac{1}{\ell}{\ell(\ell+1)-2\choose\ell-1}<0,
$$
since
$$
\frac{1}{\ell}\cdot\frac{\ell(\ell+1)-2}{\ell-1}>1
$$
and the rest of the pairwise quotients in the binomial coefficient
only make this term larger.  Thus, by the Bimodal Theorem, the desired
inequality for $n_{\ell+1}$ also holds.
\end{proof}

The next lemma will permit us to get preliminary bounds on $k_n$ which
will be needed to get the $M_n$ bounds later in this section.  

\begin{lemma}\label{kforl}
For each $k\ge 1$, the number of values of $n>\ell$ for which $k_n=k$ is
at least $\ell$ and at most $\ell+1$.
\end{lemma}
\begin{proof}
We may assume $k\ge\ell +1$ since smaller values have already been 
examined in the previous lemma.  

We begin by showing that there are at 
least $\ell$ such values of $n$.  Let $n=n_k$.  So
$\nolinebreak{c_{n-1,k}-c_{n-1,k-1}<0}$.  Since $k_n\le 
k_{n+1}$, we need only to establish that 
$\nolinebreak{c_{n+\ell-1,k+1}-c_{n+\ell-1,k}<0}$.  Hence it suffices to show that
$$
c_{n+\ell-1,k+1}-c_{n+\ell-1,k}\le c_{n-1,k}-c_{n-1,k-1}.
$$
Using~\ree{cnidiff} and Lemma~\ref{DD} with $n=k-1$, we see that the
previous inequality will follow if we can show that
\begin{eqnarray*}
\lefteqn{{k-1\choose \ell-1}-\frac{1}{\ell}{n+\ell-k-2\choose\ell-1}}
\\[5pt] 
& &\le
\frac{(\ell+1)k-1-(n-1)}{\ell}
\left[{n+\ell-k+2\choose\ell-1}-{n-k-1\choose\ell-1}\right].
\end{eqnarray*}

The case $\ell=2$ follows from straightforward computation, so we may 
assume $\ell\ge 3$ for the rest of this part of the proof.  From 
Proposition~\ref{nk-ub},
$$
\frac{(\ell+1)k-1-(n-1)}{\ell}\ge\frac{1}{\ell}>0.
$$
Because of this and the fact that $n+\ell-k+2>n-k-1$, 
the last inequality in the previous paragraph will follow if we can show
$$
{k-1\choose \ell-1}\le\frac{1}{\ell}{n+\ell-k-2\choose \ell-1}.
$$
This simplifies to
$$
\ell(k-1)\dots(k-\ell+1)\le(n+\ell-k-2)\dots(n-k).
$$
For $\ell\ge 3$ there are sufficiently many factors on both sides of
this inequality so that it will be proved if both
\beq\label{kforl-11}
k-\ell+1\le n-k
\eeq
and
\beq\label{kforl-22}
\ell(k-1)(k-2)\le(n+\ell-k-2)(n+\ell-k-3).
\eeq
Both inequalities follow from the upper bound in 
Proposition~\ref{kn-crude} as follows.  For (\ref{kforl-11}) we have 
$2k-\ell+1<2k<\ell k<n$.  For (\ref{kforl-22}), note 
$\ell(k-1)(k-2)<n(k-2)$ while
\bea
(n+\ell-k-2)(n+\ell-k-3)
&>&
(\ell k+\ell-k-2)(n+\ell-k-3)\\
&=&
((k+1)(\ell-1)-1)(n+\ell-k-3)\\
&\ge&
(2k-4)(n-k)\\
&>&
2(k-2)(n/2),
\eea
since $k<n/\ell<n/2$.  This completes the proof that for all $k$ there 
are at least $\ell$ values of $n$ for which $k_n=k$.

We would now like to show that for all $k$ there are at most $\ell+1$ 
values of $n$ for which $k_n=k$.  We do this by induction on $k$.  
Lemma~\ref{basecases} gives the result for $k\le\ell$, so we may 
assume that $k>\ell$ and that the result is true for all values less 
than $k$.  Let $n=n_k$, so $c_{n,k}-c_{n,k-1}\ge 0$.  Then by the 
Continuity Theorem it suffices to show
$$
c_{n+\ell+1,k+1}-c_{n+\ell+1,k}\ge 0
$$
since that will imply that $k_{n+\ell+1}\ge k+1$.  So it will be 
sufficient to show
$$
c_{n+\ell+1,k+1}-c_{n+\ell+1,k}\ge c_{n,k}-c_{n,k-1}.
$$
Using~\ree{cnidiff} and rearranging terms gives the equivalent 
inequality
$$
(M_{k+1}-M_k)-(M_k-M_{k-1})\ge\frac{n-(\ell+1)k+1}{\ell}
	\left[{n+\ell-k\choose\ell-1}+{n-k\choose\ell-1}\right].
$$
So by Lemma~\ref{DD}, it suffices to show
$$
n-(\ell+1)k+1=n_k-(\ell+1)k+1\le 0
$$
and this is true by Proposition~\ref{nk-ub}.
\end{proof}

Combining this result and
Proposition~\ref{nk-ub} immediately gives an upper bound  for
$n'_k$ which is defined as the 
{\it largest\/} value of $n$ such that $k_n=k$.  This will
be important for our lower bound on $k_n$ in the next section.

\bco
\label{n'k_ub}
For all $k\ge1$ we have
$$
n'_k\le(\ell+1)k+\ell-1.\qedthm
$$
\eco

We will obtain better bounds on $k_n$ in the next section by using our
upcoming bounds on the difference $M_n-M_{n-1}$.   But for the proof of
the latter result we need a weaker upper bound which comes from Lemma~\ref{kforl}.
\begin{lemma}\label{kn-first-bound}
For all $n\ge (\ell+1)\ell$ we have
$$
k_n\le\frac{n-\ell}{\ell}.
$$
\end{lemma}
\begin{proof}
Using Lemma~\ref{basecases}, it is easy to see that this result holds for
$(\ell+1)\ell\le n \le (\ell+1)^2-1$.  To finish the demonstration, it
suffices to prove the result for each $n_k$ where $k>\ell+1$.  
We may assume, by 
induction on $k$, that $k-1\le(n_{k-1}-\ell)/\ell$.  Now, by Lemma~\ref{kforl}, 
$n_k\ge n_{k-1}+\ell$ which combines with the previous inequality to complete the proof.
\end{proof}

We will also need a technical corollary of the previous lemma.
\bco
\label{technical}
For all $n\ge (\ell+1)\ell$ we have, with $k=k_n$,
$$
\beta\frac{(k-\ell+1)^\ell}{\ell!}+{n-k-1\choose \ell}
	\ge \beta\frac{k^\ell}{\ell!}+\frac{(n-k-\ell)^\ell}{\ell!}.
$$
\eco
\begin{proof}
Rearranging terms and multiplying by $\ell!$, it suffices to show
\beq
\label{sts}
(n-k-1)(n-k-2)\cdots(n-k-\ell)-(n-k-\ell)^\ell\ge
	\be\left[k^\ell-(k-\ell+1)^\ell\right].
\eeq
Now using terminating approximations for positive and alternating
series we have
$$
(n-k-1)(n-k-2)\cdots(n-k-\ell)
	\ge (n-k-\ell)^\ell +{\ell\choose2} (n-k-\ell)^{\ell-1},
$$
and
$$
(k-\ell+1)^\ell
	\ge k^\ell - \ell(\ell-1) k^{\ell-1},
$$
respectively.  Comparing these with~\ree{sts} reduces us to proving
$$
(n-k-\ell)^{\ell-1}/2\ge\be k^{\ell-1}.
$$
Lemma~\ref{whatisbeta2} gives us $\be\le 1/2$ so we will be done if
$n-k-\ell\ge k$.  But by the previous Lemma, 
$k\le (n-\ell)/\ell\le (n-\ell)/2$ which is equivalent.
\end{proof}

We are now ready to prove one of our most useful results which gives
bounds on the differences $M_n-M_{n-1}$.  This will be used to get
both our bounds 
on $M_n$ in this section and our bounds on $k_n$ in the next.

\begin{theorem}
\label{Mn-Mn-1}
If $\ell\ge3$ and $n\ge1$, then we have
$$
M_n-M_{n-1}\le \beta\frac{(n-1)^\ell}{\ell!}.
$$
Furthermore, for all $\ell\ge 2$ and $n\ge\ell$,
$$
M_n-M_{n-1}\ge \beta\frac{(n-\ell)^\ell}{\ell!}.
$$
\end{theorem}
\begin{proof}
We begin by proving the upper bound by induction on $n$.  For $n\le\ell$ 
this bound is trivial, so we may assume that $n\ge \ell+1$. 
So $k=k_n$ is well-defined.
Directly from the definitions
$$
M_{n-1}\ge c_{n-1,k}=M_k+k{n-k-1\choose\ell}.
$$
Combining this with~\ree{Mn} yields
\beq
\label{same1}
M_n-M_{n-1}\le k{n-k-1\choose \ell-1}\le \frac{k(n-k-1)^{\ell-1}}{(\ell-1)!}.
\eeq
Similarly,
$$
M_{n-1}\ge c_{n-1,k-1}=M_{k-1}+(k-1){n-k\choose\ell}.
$$
Also  $(n-k)(n-k-2)\le(n-k-1)^2$, and because $\ell\ge3$ there are
enough factors in the binomial coefficient so that
\beq\label{same2}
M_n-M_{n-1}\le M_k-M_{k-1}+{n-k\choose\ell}
\le M_k-M_{k-1}+\frac{(n-k-1)^\ell}{\ell!}.
\eeq

Combining (\ref{same1}) and (\ref{same2}) we get
\beq\label{both-cases}
M_n-M_{n-1}\le 
\gamma\left(M_k-M_{k-1}+\frac{(n-k-1)^\ell}{\ell!}\right)
+(1-\gamma)\frac{k(n-k-1)^{\ell-1}}{(\ell-1)!}
\eeq
for all  $\gamma\in[0,1]$.  
By induction $M_k-M_{k-1}\le \beta (k-1)^\ell/\ell!<\beta k^\ell/\ell!$.  
Making this substitution and setting $\gamma = k/(n-1)$ gives 
\bea
M_n-M_{n-1}
&\le&
\gamma\left(\frac{\beta\gamma^\ell (n-1)^\ell}{\ell!}
		+\frac{(1-\gamma)^\ell (n-1)^\ell}{\ell!}\right)
+(1-\gamma)\frac{\gamma(1-\gamma)^{\ell-1}(n-1)^\ell}{(\ell-1)!}\\
&=&
\left(\beta\gamma^{\ell+1}+(\ell+1)\gamma
(1-\gamma)^\ell\right)\frac{(n-1)^\ell}{\ell!}.
\eea
By Lemma~\ref{whatisbeta}, we know that for all $\gamma\in[0,1]$,
$$
\frac{(\ell+1)\gamma
(1-\gamma)^{\ell-1}}{1+\gamma+\gamma^2+\dots+\gamma^\ell}\le\beta,
$$
so
$$
(\ell+1)\gamma(1-\gamma)^\ell\le \beta(1-\gamma^{\ell+1}).
$$
It follows that
$$
\beta\gamma^{\ell+1}+(\ell+1)\gamma(1-\gamma)^\ell\le\beta,
$$
completing the proof of the upper bound.

We will have to break the proof of the lower bound into two cases
depending on the size of $n$. 

First suppose that $n<(\ell+1)\ell$.  By the Continuity Theorem we
have two subcases depending upon whether $k_{n-1}=k-1$ or $k$.
Suppose that the former is true so that we have, by
Lemma~\ref{basecases}, $k\le\ell$.
Then using $M_k=M_{k-1}=0$ and~\ree{Mn} gives
$$
M_n-M_{n-1}=k{n-k\choose\ell}-(k-1){n-k\choose\ell}={n-k\choose\ell}.
$$
We would like to show that the right-hand side of this inequality is at least 
$\beta(n-\ell)^\ell/\ell!$.  So by Lemma~\ref{whatisbeta2}, it
suffices to show that 
$$
\left(1-\frac{1}{\ell}\right)^{\ell-1}(n-\ell)^{\ell}\leq
(n-k)(n-k-1)\cdots(n-k-l+1).
$$
Note that since $k\le\ell$ we have $n-\ell\le n-k$ so we are
reduced to proving 
$$
\left(1-\frac{1}{\ell}\right)^{\ell-1}(n-\ell)^{\ell-1}\leq
(n-k-1)\cdots(n-k-l+1).
$$
This last inequality will follow if we can show
\beq
\label{bothcases}
\left(1-\frac{1}{\ell}\right)(n-\ell)\leq n-k-l+1.
\eeq
But multiplying out the left-hand side and cancelling shows that this
is true because of Proposition~\ref{kn-crude}.

Now suppose that $k_{n-1}=k$.  Then Lemma~\ref{basecases} implies that
$k<\ell$ because of the bounds on $n$ in this case.
As before, we can compute
\beq
\label{samek}
M_n-M_{n-1}=k{n-k\choose\ell}-k{n-k-1\choose\ell}=k{n-k-1\choose\ell-1}.
\eeq
Using Lemma~\ref{whatisbeta2} again, we see that we need to prove
$$
\left(1-\frac{1}{\ell}\right)^{\ell-1}(n-\ell)^{\ell}\leq
\ell k (n-k-1)(n-k-1)\cdots(n-k-l+1).
$$
But by Proposition~\ref{kn-crude} again
$$
\ell^2 k > (\ell^2-1) k =(\ell-1)(\ell+1) k \ge (\ell-1) (n-\ell)=
\ell\left(1-\frac{1}{\ell}\right)(n-\ell).
$$
Furthermore, $k<\ell$ implies $n-k-1\ge n-\ell$ and~\ree{bothcases} takes
care of the remaining factors.

We may now assume that $n\ge (\ell+1)\ell$.  Again we have two subcases.  If $k_{n-1}=k$, then $k\ge\ell$.  Also~\ree{samek} still holds and so
\beq
\label{firstdiff}
M_n-M_{n-1}\ge \frac{k(n-k-\ell)^{\ell-1}}{(\ell-1)!}
=\ell\gamma(1-\gamma)^{\ell-1}\frac{(n-\ell)^\ell}
{\ell!},
\eeq
where $\ga$ is defined by $k=\gamma(n-\ell)$.  Note that $\ga<1/\ell$ by
Lemma~\ref{kn-first-bound}.  Also, by our remarks about the function
$f(x)$ of equation~\ree{f(x)} in the proof of Lemma~\ref{whatisbeta2}, we have
the lower bound in the theorem as long as $\gamma\in [\alpha,1/\ell)$.

To see what happens if $\gamma<\alpha$, we use the fact that
$M_n\ge c_{n,k+1}$, Corollary~\ref{technical}, and induction to get
\bea
M_n-M_{n-1}
	&\ge& M_{k+1}-M_k+{n-k-1\choose \ell}\\
	&\ge& \beta\frac{(k-\ell+1)^\ell}{\ell!}+{n-k-1\choose \ell}\\
	&\ge& \beta\frac{k^\ell}{\ell!}+\frac{(n-k-\ell)^\ell}{\ell!}\\
	&=&\left(\be\ga^\ell+(1-\ga)^\ell\right)\frac{(n-\ell)^\ell}{\ell!}.
\eea
But since $\ga\in[0,\al)$, we can use
Lemma~\ref{whatisbeta} to conclude that our desired lower bound holds.
So we are now done with the case where $k_{n-1}=k$.

Now assume that $k_{n-1}=k-1$ so that $k>\ell$ because of the bound on
$n$.  Then from~\ree{Mn} we get that
$$
M_n-M_{n-1}=M_k-M_{k-1}+{n-k\choose \ell}.
$$
Thus the first bound in the previous string of inequalities holds with
$k$ replaced by $k-1$.  But $k-1\ge\ell$ so the same arguments used
there apply to give the lower bound we seek.  Similarly, since
$M_n\ge c_{n,k-1}$ we can use~\ree{c_ni} to get
$$
M_n-M_{n-1}\ge (k-1){n-k\choose \ell-1},
$$
which can be compared with~\ree{firstdiff} to complete the proof of
this case and of the theorem itself.
\end{proof}

We are now in a position to take care of most of the cases in
Theorem~\ref{MAIN}. 
\bth
\label{Mn-bounds}
Suppose $n\ge\ell$.  If $\ell\ge3$, then we have
$$
M_n\le \beta\frac{n^{\ell+1}}{(\ell+1)!}.
$$
Furthermore, for all $\ell\ge 2$,
$$
M_n\ge \beta\frac{(n-\ell)^{\ell+1}}{(\ell+1)!}.
$$
\eth
\begin{proof}
Both bounds are trivial if $n=\ell$ since $M_\ell=0$.  So suppose
$n>\ell$.

For the upper bound, we use the the previous theorem and the standard
way in which sums are used to
bound integrals to get
\bea
M_n	&=&\sum_{i=\ell+1}^n \left(M_i-M_{i-1}\right)\\
	&\le&\frac{\be}{\ell!}\sum_{i=\ell+1}^n (i-1)^\ell\\
	&\le&\frac{\be}{\ell!}\int_0^n x^\ell\ dx\\
	&=&\frac{\beta}{(\ell+1)!}n^{\ell+1}.
\eea

There are two possible proofs of the lower bound at this point.
Either one can mimic the demonstration of the upper bound or appeal to
Theorems~\ref{albert} and \ref{price} to get
$$
M_n\ge\be{n\choose \ell+1}\ge \be\frac{(n-\ell)^{\ell+1}}{(\ell+1)!}.
$$
Using either technique, we are done.
\end{proof}

\section{Bounds on $k_n$}
\label{k_n}

We can now use the results of the previous section to
supply bounds for $k_n$ which will be a considerable
improvement over those obtainable from Price's work.  The best that
can be gotten 
from Corollary~\ref{price-kn} is $k_n = \al n + o(n)$.  We will prove
that in fact $k_n = \al n + O(1)$ with a constant inside the big oh
that is less than 2.

\begin{theorem}
\label{bruce}
For $\ell\ge 3$ and $n>\ell$ we have 
$$
k_n\le\alpha(n-\ell)+1.
$$
\end{theorem}
\begin{proof}
Let $k=k_n$.  Note that it suffices to prove the bound when $n=n_k$.

Clearly the result is true for $k=1$ and $n_1=\ell+1$.
Next suppose that $2\le k\le\ell+1$.  Then by Lemma~\ref{basecases},
our desired inequality is equivalent to $k-1\le(\ell+1)\al(k-1)$ which
is true by~\ree{alpha-bounds}. 

For $k>\ell+1$ we still have Proposition~\ref{nk-ub} which gives
$n-(\ell+1)k+1\le0$.  So by Theorem~\ref{Mn-Mn-1} and the fact that
$\ell\ge 3$, 
\bea
0
&\le&
c_{n,k}-c_{n,k-1}\\
&=&
M_k-M_{k-1}+\frac{n-(\ell+1)k+1}{\ell}{n-k\choose \ell-1}\\
&\le&
\beta\frac
{(k-1)^\ell}{\ell!}+\frac{[n-\ell-(\ell+1)(k-1)][n-\ell-(k-1)]^{\ell-1}}{\ell!
}.
\eea

Define $\gamma$ by $k-1=\gamma(n-\ell)$.  So it suffices to show 
$\gamma\le\alpha$.  Clearly $\gamma\ge 0$ and since $k>\ell+1$ we
can apply Lemma~\ref{kn-first-bound} to get 
$\gamma(n-\ell)<k\le (n-\ell)/\ell$, so  
$\gamma<1/\ell$.  Rewriting the last expression in the previous
paragraph in terms of $\ga$ and cancelling $\ell!$ gives
$$
0\le\beta\gamma^\ell
(n-\ell)^\ell+(1-(\ell+1)\gamma)(1-\gamma)^{\ell-1}(n-\ell)^\ell.
$$
Since $n>\ell$ we have
$$
0\le\beta\ga^\ell+(1-(\ell+1)\gamma)(1-\gamma)^{\ell-1}.
$$
Now define
$$
g(x)=\beta x^\ell+(1-(\ell+1)x)(1-x)^{\ell-1},
$$
so we have $g(\gamma)\ge 0$.  Using the defining equations for $\al$
and $\be$ gives
$$
g(\al)=\ell\alpha^{\ell+1}
(1-\alpha)^{\ell-1}+(1-(\ell+1)\alpha)(1-\alpha)^{\ell-1}=0.
$$
This implies that
$$
0\le g(\ga)=g(\al)+\int_\al^\ga g'(x)\ dx = \int_\al^\ga g'(x)\ dx.
$$
Since $0\le\alpha,\gamma\le 1/\ell$, we can prove $\gamma\le\alpha$ by showing 
that $g'(x)<0$ on $[0,1/\ell]$.  Now
$$
g'(x)=\beta\ell 
x^{\ell-1}-(\ell+1)(1-x)^{\ell-1}-(1-(\ell+1)x)(\ell-1)(1-x)^{\ell-2}.
$$
So we want, after transposing terms,
$$
\beta\ell x^{\ell-1}
<(\ell+1)(1-x)^{\ell-1}+(1-(\ell+1)x)(\ell-1)(1-x)^{\ell-2}.
$$
Taking the maximum of the left-hand side and the minimum of the 
right-hand side on the interval $[0,1/\ell]$, it suffices to show that
$$
\frac{\ell^2\alpha(1-\alpha)^{\ell-1}}{\ell^{\ell-1}}
<
\frac{(\ell+1)(\ell-1)^{\ell-1}}{\ell^{\ell-1}}
   -\frac{(\ell-1)^{\ell-1}}{\ell^{\ell-1}}
=
\frac{\ell(\ell-1)^{\ell-1}}{\ell^{\ell-1}}.
$$
But $\ell^2\alpha(1-\alpha)^{\ell-1}<\ell(\ell-1)^{\ell-1}$ since 
by equation~\ree{alpha-bounds} we have both that $\ell^2\alpha<\ell$
and that $(1-\alpha)^{\ell-1}<(\ell-1)^{\ell-1}$. 
\end{proof}

We also have a lower bound with only a slightly larger constant.

\bth
\label{martin}
For $\ell\ge 2$ and $n$ sufficiently large, 
$$
k_n\ge\alpha(n-\ell)-1.
$$ 
\eth
\begin{proof}
Let $k=k_n$.  Note that it suffices to prove the bound when $n=n'_k$.
Note also that by Corollary~\ref{n'k_ub} we have $n-(\ell+1)k-\ell<0$. 
Using this fact, the Bimodal Theorem, equation~\ree{cnidiff}, and
Theorem~\ref{Mn-Mn-1}, we have
\bea
0&>&c_{n,k+1}-c_{n,k}\\
 &=&M_{k+1}-M_k+\frac{n-(\ell+1)k-\ell}{\ell}{n-k-1\choose\ell-1}\\
 &\ge&\frac{\be(k-\ell+1)^\ell+[n-(\ell+1)k-\ell](n-k-1)^{\ell-1}}{\ell!}.\\
\eea 

Now define $\ga$ by $k-\ell+1=\ga(n-\ell)$. By taking $n$ sufficiently
large we can assume that $k+1\ge\ell$ and so $\ga\ge0$.  
Also, Theorem~\ref{bruce} implies that
$\ga(n-\ell)\le\al(n-\ell)-\ell+2$ and so $\ga\le\al$.
Substituting for $\ga$ to replace $k$ in the last inequality of the previous
paragraph we get, after multiplying by $\ell!/(n-\ell)^\ell$,
$$
\be\ga^\ell+
\left[1-\frac{\ell^2-1}{n-\ell}-(\ell+1)\ga\right](1-\ga)^{\ell-1}
<0.
$$ 
Let $\ep=(\ell^2-1)/(n-\ell)$ and note that we can make $\ep$ as small
a positive number as we wish by taking $n$ large.  Define a function
$$
h(x)=\be x^\ell+[1-\ep-(\ell+1)x](1-x)^{\ell-1}
$$
so that $h(\ga)<0$.  Using the defining equations for $\al$ and $\be$,
one can also compute that
\beq
\label{hal}
h(\al)=-\ep(1-\al)^{\ell-1}.
\eeq

We want to mimic the integration trick used in the proof of the upper
bound for $k_n$, so we need some information about $h'(x)$.  First
note that
\bea
h'(x)&=&\ell\be x^{\ell-1}
   -\left[(\ell+1)(1-x)+(\ell-1)(1-\ep-(\ell+1)x)\right](1-x)^{\ell-2}\\
&=&\ell\be x^{\ell-1}
   -\left[(\ell+1)(1-\ell x)+(\ell-1)(1-\ep)\right](1-x)^{\ell-2}.
\eea
Taking one more derivative, one can see that $h''(x)\ge0$ on
$[0,1/\ell]$ as long as the factor in the 
final set of square brackets above is nonnegative.  And this can be
ensured by taking $1-\ep>0$.  So $h'(x)$ is increasing on this interval, and since
$\ga\le\al< 1/\ell$ we can write  
\beq
\label{int}
0>h(\ga)=h(\al)-\int_\ga^\al h'(x)\ dx \ge h(\al) - (\al-\ga)h'(\al).
\eeq

Next we claim that
$$
h'(\al)\le -(\ell-1)(1-\al)^{\ell-2}.
$$
Using the second expression for $h'(x)$ and the definition of $\be$,
we see that it is sufficient to prove, after cancelling
$(1-\al)^{\ell-2}$, that
$$
\ell^2\al^\ell(1-\al)-(\ell+1)(1-\ell\al)-(\ell-1)(1-\ep)\le -(\ell-1).
$$ 
Expanding the left-hand side and using
$\ell\al^{\ell+1}=(\ell+1)\al-1$ on the $-\ell^2\al^{\ell+1}$ term
reduces this inequality, after massive cancellation, to
\beq
\label{ep}
\ell^2\al^\ell+(\ell-1)\ep\le 1.
\eeq
But this last equation is true for sufficiently small $\ep$ since,
by~\ree{alpha-bounds}, 
$$
\ell^2\al^\ell\le\ell^2\al^2<1.
$$
So we have proved the claim.

Now divide~\ree{int} by $h(\al)$ (which is negative by~\ree{hal}) and
use the claim as well as~\ree{alpha-bounds} again to get
\bea
0&<&1-(\al-\ga)\frac{h'(\al)}{h(\al)}\\
 &\le& 1-(\al-\ga)\frac{n-\ell}{(\ell+1)(1-\al)}\\
 &<& 1-(\al-\ga)\frac{n-\ell}{\ell}.
\eea
Solving for $\ga$ in this last inequality and plugging into its
defining equation gives
\bea
k&=&\ga(n-\ell)+\ell-1\\
 &\ge&\left(\al-\frac{\ell}{n-\ell}\right)(n-\ell)+\ell-1\\
 &=&\al(n-\ell)-1
\eea
as desired.
\end{proof}

To give a feel for how good these bounds are, we prove the
following corollary.
\bco
For $\ell\ge3$ and $n>\ell$ we have
$$
k_n-\al n< 1/4.
$$
For $\ell\ge2$ and sufficiently large
$n$ we have
$$
k_n-\al n > -2.
$$
\eco
\begin{proof}
The lower bound follows immediately from the previous theorem
and~\ree{alpha-bounds}.  For the upper bound, it is easy to show by
taking second derivatives that $f_{\ell-1}(x)\ge f_\ell(x)$ on the
interval $[0,1/\ell]$.  It follows that $\al$ is a decreasing function
of $\ell$.  Furthermore, using~\ree{alpha-bounds} again shows that
$|1-\ell\al|<\al$.  Combining these observations with Theorem~\ref{bruce} gives
$$
k_n-\al n\le 1-\ell\al<\al.
$$ 
But now we are done since $1-\ell\al<1/4$ when $\ell=3$ and $\al<1/4$
for $\ell\ge4$.
\end{proof}

\section{The upper bounds for $\ell=2$}
\label{ell=2}

To complete the proof of Theorem~\ref{MAIN} we must address the upper
bound when $\ell=2$.  This result, as well as the upper bound on $k_n$
in the previous section, depends on Theorem~\ref{Mn-Mn-1} where the
restriction $\ell\ge3$ first appeared.  This is not an accident as
that theorem is false for $\ell=2$.  For example, when $\ell=2$ 
we have $M_{17} - M_{16} = 60$, but
$\beta 16^2/2$ is approximately $59.405$.  Worse yet, our computer experiments have shown that this is not an isolated counterexample.  However, a weaker upper bound is true.
\begin{theorem}
\label{Mndiff}
For $\ell=2$ and $n\ge\ell$ we have
$$
M_n-M_{n-1}\le \beta\frac{n^2}{2}.
$$
\end{theorem}
\begin{proof}
The proof is very similar to the demonstration of
Theorem~\ref{Mn-Mn-1}.  There 
are only two changes.  The first is that when bounding binomial
coefficients one uses powers of $n-k$ rather than $n-k-1$.  Note that
this removes the necessity to have $\ell\ge3$.  The other modification
is that one substitutes $\ga=k/n$.  The rest of the proof proceeds as before.
\end{proof}

We can now obtain the $\ell=2$ upper bound in Theorem~\ref{MAIN}.  One
uses the same proof as Theorem~\ref{Mn-bounds} but with the
previous result taking the place of Theorem~\ref{Mn-Mn-1}.  Because of
the similarilty, we omit the details.
\bth
For $\ell=2$ and $n\ge\ell$ we have
$$
M_n\le \beta\frac{(n+1)^3}{3!}.\qedthm
$$
\eth

To obtain the bounds on $k_n$ in this case, note that $f_\ell(x)$
always has
$x=1$ as a root.  So dividing $f_2(x)$ by $x-1$, we see that $\al$
must satisfy
\beq
\label{al}
\al^2=\frac{-2\al+1}{2}.
\eeq
We can now plug this into the defining equation for $\be$ to get
\beq
\label{be}
\be=4\al-1.
\eeq
\bth
For $\ell=2$ and all $n\ge3$ we have
$$
k_n\le \al n + 1/2.
$$
\eth
\begin{proof}
Let $k=k_n$ as usual.
Using Theorem~\ref{Mndiff}, as well as
equations~\ree{cnidiff} and~\ree{be} gives
\bea
0&\le&c_{n,k}-c_{n,k-1}\\
 &=&M_k-M_{k-1}+\frac{n-3k+1}{2}{n-k\choose1}\\
 &\le& (4\al-1)\frac{k^2}{2} + \frac{(3k-n-1)(k-n)}{2}.
\eea

Let
\bea
f(x)&=&(4\al-1)x^2 + (3x-n-1)(x-n)\\
	&=&(4\al+2)x^2-(4n+1)x+(n^2+n).
\eea
The vertex of this parabola is at $x_0=(4n+1)/(8\al+4)$ and from
Proposition~\ref{kn-crude} we have $k<n/2<x_0$.  Combining this with the
fact that  $f(k)\ge0$ shows that $k$ is at most the smaller of the two
roots of $f(x)$ which is
$$
r=\frac{4n+1-\sqrt{(4n+1)^2-4(n^2+n)(4\al+2)}}{8\al+4}.
$$
To complete the proof we need to show that $r\le \al n + 1/2$.
Rearranging terms in this last inequality and using~\ree{al} shows
that we need to prove
$$
\sqrt{(4n+1)^2-4(n^2+n)(4\al+2)}\ge 4\al n - 4\al-1.
$$ 
Since $n\ge3$, the right-hand side of this last inequality is
positive.  So we can square it and use~\ree{al} again to reduce our
task to proving $(16-40\al)n+(8\al-8)\ge0$.  But this is true since
$n\ge3$ and the theorem is proved.
\end{proof}

\section{Comments and conjectures}

There are several ways in which this work could be continued.  We list
some of them here in the hopes that the reader will be interested.

\medskip

1.  We have already noted that  the upper bound in
Theorem~\ref{Mn-Mn-1} is not true for $\ell=2$.  However, numerical
evidence indicates that the succeeding results are still valid, even though
the proofs we have given will not work.  In particular, we make the
following conjecture.
\bcon
For $\ell=2$ and $n>\ell$ we have
$$
M_n\le\be \frac{n^{\ell+1}}{(\ell+1)!}
$$
and
$$
k_n \le \al (n-\ell) + 1.
$$
\econ

\medskip

2.  The lower bound given for $k_n$ in Theorem~\ref{martin} suffers
from the fact that our demonstration only works for sufficiently large
$n$.  The most restrictive place where this is used is in the proof
that inequality~\ree{ep} holds and there we need $n$ to be at least on
the order of $\ell^3$.  But numerical calculations suggest that an
even better bound holds for all $n$.
\bcon
For all $\ell\ge2$ and $n>\ell$ we have
$$
k_n\ge\al(n-\ell).
$$
\econ
We should note that there are examples where $k_n$ is not the closest
integer to $\al n$.  So, given the upper bound we have already proven, one can not hope to substantially improve upon this conjecture.

\medskip

3.  The reader will have noticed that the Continuity Theorem has been of
fundamental importance in proving the results in this paper.  This
leads us to wonder if something can be said for a larger class of
layered patterns $q$.  By Theorem~\ref{strom-layered}, one can
still define $k_n$ as the maximum 
length of the word remaining after removing the last layer of a
$q$-optimal layered $n$-permutation.  So we would like to be able to
say something about the sequence $k_n$.  There are some results in this
regard in Price's thesis~\cite{price} for patterns with at most two
layers and certain patterns with all layer lengths two.

Another of our main tools which might be amenable to generalization to
other layered permutations is the Bimodal Theorem.  One can still
define $c_{n,i,q}$ to be the maximum number of copies of $q$  in a
layered $n$-permutation where the last layer has length $n-i$.
Knowing the shape of the sequence $(c_{n,i,q})_{0\le i< n}$ could be
useful in getting information about the packing density of $q$.

\medskip

4.  Because of Theorem~\ref{albert}, it is easy to generalize the
lower bound of Theorem~\ref{MAIN} to all patterns.  The proof is the
same as the second proof of the lower bound in Theorem~\ref{Mn-bounds}
and so is left to the reader.
\bth
If $q$ is a pattern of length $L$ and $n\ge L$ then
$$
M_{n,q}\ge \de(q)\frac{(n-L+1)^L}{L!}.\qedthm
$$
\eth
We conjecture that the corresponding upper bound holds as well.
\bcon
If $q$ is a pattern of length $L$ and $n\ge L$ then
$$
M_{n,q}\le \de(q)\frac{n^L}{L!}.
$$
\econ

\medskip

5.  Finally, we should point out that since Herb Wilf first defined packing densities in 1992 at the SIAM meeting on Discrete Mathematics, only packing densities of layered permutations (or permutations equivalent to layered permutations under one of the 8 routine symmetries) have been computed.  The first open cases are of length four, where Albert, Atkinson, Handley, Holton, and Stromquist~\cite{aahhs:opdop} gave the bounds
$$
0.19657\le\delta(1342)\le 2/9,
$$
and
$$
51/511\le\delta(2413)\le 2/9.
$$
While we are hopeful that the approach presented in this paper (and in particular, generalizations of the Continuity and Bimodal Theorems) may prove fruitful in other layered cases, our approach seems to offer no additional hope in the nonlayered cases.

\begin{\bib}{99}

\bibitem{aahhs:opdop} M. Albert, M. Atkinson, C. Handley, D. Holton, and W. Stromquist, On packing densities of permutations, {\it Electronic J.\ Combin.\/} {\bf 9} (2002), \#R5.

\bibitem{bon:iz} M. B\'ona, B. Sagan, and V. Vatter, Pattern 
frequency sequences and internal zeros, {\it Adv.\ Appl.\ Math.\/} {\bf 
28} (2002), 395-420.

\bibitem{bhm:words} A. Burstein, P. H\"{a}st\"{o}, T. Mansour, Packing patterns into words, {\it Electronic J. Combin.\/} {\bf 9(2)} (2003), \#R20.

\bibitem{h:pd} P. H\"{a}st\"{o}, The packing density of other layered 
permutations, {\it Electronic J. Combin.\/} {\bf 9(2)} (2002), \#R1.

\bibitem{price} A. Price, ``Packing densities of layered patterns,'' 
Ph.D. thesis, University of Pennsylvania, Philadelphia, PA, 1997.

\bibitem{str:plp} W. Stromquist, Packing layered posets into posets, 
manuscript.

\end{\bib}

\end{document}